\documentclass{article}

{}
{}
{}
\usepackage{color}
\usepackage{multirow}
\usepackage{makecell}
\usepackage{float}
\usepackage{setspace}
\usepackage{array}
\usepackage{longtable}
\usepackage{supertabular}
\usepackage{geometry}
\usepackage{fancyhdr} 
\usepackage{mathtools}
\usepackage{amssymb}
\usepackage{graphicx}
\usepackage{lineno,hyperref}

\usepackage{amsfonts}
\usepackage{comment}

\begin{document}

%\supertitle{Research Article}

\title{Investigating Continuous Power Flow Solutions of IEEE-14 Bus System}

\author{Bin Liu\thanks{School of Electrical Engineering and Telecommunications, The University of New South Wales, Sydney 2052, Australia and previously with Jibei Electric Power Dispatching and Control Center, State Grid Cooperation of China, Beijing 100053, China. Email address: eeliubin@hotmail.com.}, Feng Liu\thanks{Department of Electrical Engineering, Tsinghua University, Beijing 100084, China.}, Bingxu Zhai\thanks{Jibei Electric Power Dispatching and Control Center, State Grid Corporation of China, Beijing 100054, China.}, Haibo Lan\thanks{Jibei Electric Power Dispatching and Control Center, State Grid Corporation of China, Beijing 100054, China.}}
%\thanks{Department of Electrical Engineering, Tsinghua University, Beijing 100084, China}
%\author{\au{Bin Liu$^{1,2\corr}$}, \au{Wei Wei$^{2}$}, \au{Feng Liu$^{2}$}}
%\add{2}{Department of Electrical Engineering, Tsinghua University, Beijing 100084, China}
%\add{3}{State Grid Jibei Electric Power Dispatching and Control Center, Beijing 100053, China}
%\email{eeliubin@hotmail.com}}

\maketitle
\thispagestyle{fancy}          %更改plain状态
\fancyhead{}                      %清除以前的命令
\lhead{This paper is a postprint of a paper accepted by IEEJ Transactions of Electrical and Electronics Engineering and is subject to The Institute of Electrical Engineers of Japan. The copy of record is available at https://onlinelibrary.wiley.com/doi/abs/10.1002/tee.22773}           %左上角添加
\chead{}
\rhead{}
\lfoot{}
\cfoot{}   %current page number
\rfoot{\thepage}
\renewcommand{\headrulewidth}{0pt}       %把页眉线的宽度设为零，即去掉页眉线
\renewcommand{\footrulewidth}{0pt}

\begin{abstract}
This letter focuses on the multiplicity of power flow (PF) equations and presents two continuous solutions for widely studied IEEE-14 bus system. The continuous solutions are located by a method combining the semidefinite program (SDP) relaxation and reformulation linearization technique (RLT). Although the observation is non-trivial, it is of interest to researchers investigating the geometry or multiplicity nature of PF equations.
\end{abstract}

\section{Introduction}
Solving PF equations is one of the most fundamental problem in power system and it is widely reported that multiple PF solutions could exist due to the nonlinearity \cite{ref-1}-\cite{ref-5}. Recent reported methods to locate multiple PF solutions include continuous power flow-based algorithm \cite{ref-2}, the continuation method \cite{ref-3}, the homotopy continuation method \cite{ref-4} and convex relaxation based method \cite{ref-5}. In \cite{ref-4}, 30 {\it isolated} solutions are successfully located for the IEEE-14 bus system. Apart from the reported {\it isolated} PF solutions, one interesting question is whether a {\it continuous} solution exists with the slack bus fixed. If any continuous PF solution exists, it implies that there will be innumerous PF solutions, or that at least a PF solution curve can be obtained, from the mathematical perspective, which is of great interest to the community investigating geometry or multiplicity nature of PF equations.  Based on the method presented in \cite{ref-5}, this letter investigates the PF solutions of IEEE-14 bus system and finds that at least two continuous PF solutions exist. 

%Section II
\section{Formulating Power Flow Equations}
Typically, the PF problem can be formulated as a set of equations as follows, where the voltage magnitude $|V|$ and voltage angle $\theta$ of $PQ$ and $PV$ buses are state variables and reactive power generations of $PV$ and $V\theta$ buses, active power generation of $V\theta$ bus are decision variables. 
\begin{eqnarray}
\label{eq-1-1}
P_{m}^g-P_m^d=|V_{m}|^{2}g^s_{m}+\sum_{n\in \mathcal{B}_m}{P_{mn}},\quad \forall m \in \mathcal{PV \cup PQ}\\
\label{eq-1-2}
Q_{m}^g-Q_m^d=|V_{m}|^{2}b^s_{m}+\sum_{n\in \mathcal{B}_m}{Q_{mn}},\quad \forall m \in \mathcal{PQ}
\end{eqnarray}
with
\begin{eqnarray} 
\label{eq-2-1}
P_{mn}=(|V_m|^2-|V_{m}||V_{n}|\cos{\theta_{mn}})g_{mn}-|V_{m}||V_{n}|b_{mn}\sin{\theta_{mn}}\nonumber\\
\label{eq-2-2}
Q_{mn}=-|V_{m}||V_{n}|g_{mn}\sin{\theta_{mn}}-(|V_m|^2-|V_{m}||V_{n}|\cos{\theta_{mn}})b_{mn}\nonumber
\end{eqnarray}
where $P_{mn}$, $Q_{mn}$ are active and reactive power flows of line $mn$ (from bus $m$ to bus $n$); $g_{mn}$, $b_{mn}$ are conductance and susceptance of line $mn$; $P_m^g$, $Q_m^g$, $P_m^d$ and $Q_m^d$ are active and reactive power generations, active and reactive loads of bus $m$; $g^s_{m}$ and $b^s_{m}$ are grounding conductance and susceptance of bus $m$; $\mathcal{PV}$, $\mathcal{PQ}$ are the sets of $PV$ and $PQ$ buses and $\mathcal{B}_m$ is the set of buses connecting to bus $n$ excluding bus $n$ itself.

\section{Solve PF by a Relaxation Method}
Using $e_m=|V_m|\cos{\theta_m}$ and $f_m=|V_m|\sin{\theta_m}$, PF problem can be reformulated in rectangular coordinate as a quadratic constrained program (QCP) problem, whose objective is to minimize the sum of positive slack variables $s_k^+$ and $s_k^-$. QCP for PF problem can be expressed in compact form as (referred as QCPF)
\begin{eqnarray}
\label{eq-3-1}~~~~~~~~~~S_{opt} =\min \sum{(s_k^++s_k^-)}\\
\label{eq-3-2}~~~~~~~~~~\mbox{tr}(XZ_k )+s_{k}^+ -s_{k}^- =c_k,~~\forall k\\
\label{eq-3-3}~~~~~~~~~~s_{k}^+ ,s_{k}^- \ge 0,\  \forall k  \\
\label{eq-3-4}~~~~~~~~~~x^l\le x =\left[\makecell{e\\f}\right]\le x^u\\
\label{eq-3-5}~~~~~~~~~~X = x x^T
\end{eqnarray}
where ${\mbox {tr}} (\cdot)$ is the matrix trace operator; $Z_k,c_k$ are known constants related to power system parameters; $x^l/x^u$ is estimated lower/upper bound of $x$ and $\left[\makecell{e;f}\right]=[e_1,\cdots,e_N,f_1,\cdots,f_N]^T$ with $N$ denoting the number of bus in the systemystem.

QCPF is a non-convex problem due to \eqref{eq-3-5} and each optimal solution of QCPF with a zero objective value corresponds to a solution of original PF problem. Based on the reformulation linearization technique (RLT) and the SDP relaxation technique, QCPF is relaxed to a SDP problem (referred as $\text{QCPF}_{\text{sdp}}$) given below \cite{ref-5}.
\begin{eqnarray}
\label{eq-4-1}S_{cvx}=\min{\{\eqref{eq-3-1}|\eqref{eq-3-2}-\eqref{eq-3-3}, X\succeq 0, (x,X)\in \mbox{RLT}(x^l,x^u)\}}
\end{eqnarray}
where
\begin{equation}\label{eq-4-2}\begin{split}
\mbox {RLT}(x^l,x^u)=\left\{\makecell{x^l\le x \le x^u\\X\succeq0} \left| 
\begin{array}{c}
{x^l x^T + x (x^l)^T - x^l (x^l)^T \le X}  \\
{x^u x^T + x (x^u)^T - x^u (x^u)^T \le X}  \\
{X \le x (x^u)^T + x^l x^T- x^l (x^u)^T}  \\
\end{array} \right. \right\}
\end{split}\end{equation}
and $X\succeq 0$ means matrix $X$ is positive semi-definite.

By successively narrowing upper and lower bounds of $x$, the relaxation is gradually tightened and all optimal solutions with zero objective value for QCPF can be found  by solving a series of resulted SDP problems. Details of the method can be found in \cite{ref-5} and are omitted here for simplicity.

%Section II
\section{Investigating IEEE-14 Bus System}
\begin{figure}[ht]
	\centering
	\includegraphics[scale=0.25]{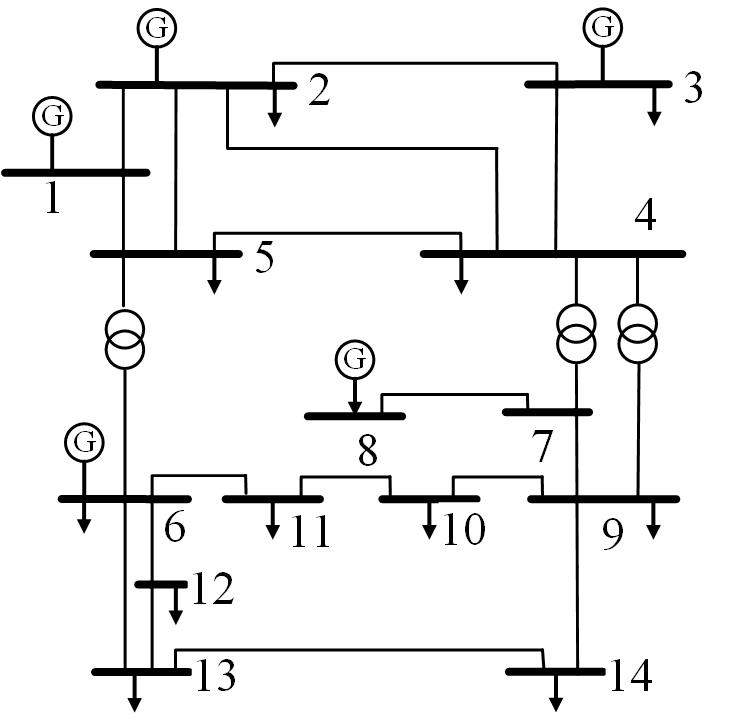}
	\caption{Topology of the IEEE 14-bus system}
	\label{fig-1}
\end{figure}

The topology and generation parameters of IEEE 14-bus system are presented in Fig.\ref{fig-1} and Table \ref{tab-2}. Other parameters can be found in \cite{ref-4} or Matpower \cite{ref-6}. Through investigating the topology of the system, an interesting characteristic of the IEEE 14-bus system can be observed, as we explain next. 
\begin{table}[ht]\renewcommand\arraystretch{0.7}
	\caption{Bus Parameters of IEEE 14-Bus System}
	%\vspace{1pt}
	\centering
	\begin{tabular}{>{\scriptsize}c>{\scriptsize}c>{\scriptsize}c>{\scriptsize}c>{\scriptsize}c>{\scriptsize}c>{\scriptsize}c>{\scriptsize}c}
		\hline
		Bus	&Type&\makecell{$P^d$ \\(MW)}&\makecell{$Q^d$\\(MVar)}&\makecell{$|V|$\\(p.u.)}&\makecell{$\theta$\\$(^\text{o})$}&\makecell{$P^g$\\(MW)}&\makecell{Shunt\\(p.u)}	\\
		\hline
		1	&	$V\theta$	&	0.0	&	0.0	&	1.060	&	0.0	&	$-$	&	$-$	\\
		2	&	$PV$	&	21.7	&	12.7	&	1.045	&	$-$	&	232.4	&	$-$	\\
		3	&	$PV$	&	94.2	&	19.0	&	1.010	&	$-$	&	40.0	&	$-$	\\
		4	&	$PQ$	&	47.8	&	-3.9	&	$-$	&	$-$	&	$-$	&	$-$	\\
		5	&	$PQ$	&	7.6	&	1.6	&	$-$	&	$-$	&	$-$	&	$-$	\\
		6	&	$PV$	&	11.2	&	7.5	&	1.070	&	$-$	&	0.0	&	$-$	\\
		7	&	$PQ$	&	0.0	&	0.0	&	$-$	&	$-$	&	$-$	&	$-$	\\
		8	&	$PV$	&	0.0	&	0.0	&	1.090	&	$-$	&	0.0	&	$-$	\\
		9	&	$PQ$	&	29.5	&	16.6	&	$-$	&	$-$	&	$-$	&	$j$0.19	\\
		10	&	$PQ$	&	9.0	&	5.8	&	$-$	&	$-$	&	$-$	&	$-$	\\
		11	&	$PQ$	&	3.5	&	1.8	&	$-$	&	$-$	&	$-$	&	$-$	\\
		12	&	$PQ$	&	6.1	&	1.6	&	$-$	&	$-$	&	$-$	&	$-$	\\
		13	&	$PQ$	&	13.5	&	5.8	&	$-$	&	$-$	&	$-$	&	$-$	\\
		14	&	$PQ$	&	14.9	&	5.0	&	$-$	&	$-$	&	$-$	&	$-$	\\
		\hline
	\end{tabular}
	\centering
	\label{tab-2}
\end{table}
%\begin{table}[ht]\renewcommand\arraystretch{0.7}
%	\caption{Line Parameters of IEEE 14-Bus System}
%	%\vspace{6pt}
%	\centering
%	\begin{tabular}{>{\scriptsize}c>{\scriptsize}c>{\scriptsize}c>{\scriptsize}c>{\scriptsize}c}
%		\hline
%		Line   &From bus  &  	To bus	&  \makecell{Impedance\\(p.u.)}	&  \makecell{Shunt\\(p.u.)} \\
%		\hline
%		1	&	1	&	2	&	0.01938 	+	$j$	0.05917 	&	0.0528 	\\
%		2	&	1	&	5	&	0.05403 	+	$j$	0.22304 	&	0.0492 	\\
%		3	&	2	&	3	&	0.04699 	+	$j$	0.19797 	&	0.0438 	\\
%		4	&	2	&	4	&	0.05811 	+	$j$	0.17632 	&	0.0340 	\\
%		5	&	2	&	5	&	0.05695 	+	$j$	0.17388 	&	0.0346 	\\
%		6	&	3	&	4	&	0.06701 	+	$j$	0.17103 	&	0.0128 	\\
%		7	&	4	&	5	&	0.01335 	+	$j$	0.04211 	&	0.0000 	\\
%		8	&	4	&	7	&	0.00000 	+	$j$	0.20912 	&	0.0000 	\\
%		9	&	4	&	9	&	0.00000 	+	$j$	0.55618 	&	0.0000 	\\
%		10	&	5	&	6	&	0.00000 	+	$j$	0.25202 	&	0.0000 	\\
%		11	&	6	&	11	&	0.09498 	+	$j$	0.19890 	&	0.0000 	\\
%		12	&	6	&	12	&	0.12291 	+	$j$	0.25581 	&	0.0000 	\\
%		13	&	6	&	13	&	0.06615 	+	$j$	0.13027 	&	0.0000 	\\
%		14	&	7	&	8	&	0.00000 	+	$j$	0.17615 	&	0.0000 	\\
%		15	&	7	&	9	&	0.00000 	+	$j$	0.11001 	&	0.0000 	\\
%		16	&	9	&	10	&	0.03181 	+	$j$	0.08450 	&	0.0000 	\\
%		17	&	9	&	14	&	0.12711 	+	$j$	0.27038 	&	0.0000 	\\
%		18	&	10	&	11	&	0.08205 	+	$j$	0.19207 	&	0.0000 	\\
%		19	&	12	&	13	&	0.22092 	+	$j$	0.19988 	&	0.0000 	\\
%		20	&	13	&	14	&	0.17093 	+	$j$	0.34802 	&	0.0000 	\\
%		\hline
%	\end{tabular}
%	\centering
%	\label{tab-1}
%\end{table}

Obviously, neither active/reactive load nor grounding admittance is connected to bus 7 and bus 8, and  no shunt impedance is connected to the line connecting the two buses. As no active power is generated from  bus 8, the following equations can be derived according to \eqref{eq-1-1}-\eqref{eq-1-2}.
\begin{eqnarray}
\label{eq-5-1}
P^g_{8}=P_{87}=-|V_{8}||V_{7}|b_{87}\sin({\theta_{8}-\theta_{7}})=0\\
\label{eq-5-2}
Q^g_{8}=Q_{87}=-|V_{8}|^{2}b_{87}+|V_{8}||V_{7}|b_{87}\cos({\theta_{8}-\theta_{7}})
\end{eqnarray}

At least one of the following equations holds according to \eqref{eq-5-1} due to $|V_{8}|>0$.
\begin{eqnarray}
\label{eq-6-1}
~~~~~~~~~~~~~~~~~~~~~~~~~|V_7|=0\\
\label{eq-6-2}
~~~~~~~~~~~~~~~~~~~~~~~~~\sin{(\theta_{8}-\theta_{7})}=0
\end{eqnarray}

We focus on \eqref{eq-6-1}, which means bus 7 is a grounding point partitioning the whole system into two subsystems. We here denote the subsystem containing bus 7, bus 8 and the line connecting them as $S_1$ and the other subsystem as $S_2$. With simple calculation, it can be easily obtained that $Q^g_{8}=-|V_{8}|^{2}b_{87}=6.7448~ \text{p.u.}$. Further, we have 
\begin{eqnarray}
\label{eq-7-1}
e_7=|V_7|\cos{\theta_7}=0;f_7=|V_7|\sin{\theta_7}=0\\
\label{eq-7-3}
e_8=|V_8|\cos{\theta_8}=1.06\cos{\theta_8}\\
\label{eq-7-4}
f_8=|V_8|\sin{\theta_8}=1.06\sin{\theta_8}
\end{eqnarray}

Therefore, arbitrary value of $\theta_8$ can always satisfy the PF equations of $S_1$, which means that any PF solution of $S_2$, if any, will be a continuous PF solution of the IEEE-14 bus system. In other words, one or more PF solution curves can be obtained and its number of PF solution will be innumerous. With the SDP relaxation based method, two PF solutions are found for $S_2$, which are presented in Table \ref{tab-3}.

\begin{table}[ht]\renewcommand\arraystretch{0.7}
	\caption{Power Flow Solutions ($|V_7|=0$)}
	%\vspace{6pt}
	\centering
	\begin{tabular}{>{\scriptsize}c>{\scriptsize}c>{\scriptsize}c>{\scriptsize}c>{\scriptsize}c}
		\hline
		\multirow{2}*{Bus} & \multicolumn{2}{c}{$|V|$(p.u.)} & \multicolumn{2}{c}{$\theta(^\text{o})$}\\
		&$1^{\text{st}}$ Solution   & $2^{\text{nd}}$ Solution& $1^{\text{st}}$	Solution& $2^{\text{nd}}$ Solution\\
		\hline
		1	&	1.0600 	&	1.0600	&	0.0000	&	0.0000	\\
		2	&	1.0450	&	1.0450	&	-9.6493	&	-8.4330	\\
		3	&	1.0100	&	1.0100	&	-21.7179	&	-19.6333	\\
		4	&	0.7270	&	0.7525	&	-14.7476	&	-12.5904	\\
		5	&	0.7998	&	0.8309	&	-16.7431	&	-14.0091	\\
		6	&	1.0700	&	1.0700	&	-50.5993	&	-39.2774	\\
		7	&	0.0000	&	0.0000	&	$-$		&	$-$		\\
		8	&	1.0600	 &	1.0600 	&	$-$ 	&	$-$	\\
		9	&	0.1090	&	0.2592	&	-70.2622	&	-37.2568	\\
		10	&	0.2420	&	0.3814	&	-59.5697	&	-39.0664	\\
		11	&	0.6378	&	0.7118	&	-51.9195	&	-38.9400	\\
		12	&	0.9824	&	0.9950	&	-52.4829	&	-40.7917	\\
		13	&	0.9008	&	0.9280	&	-51.9006	&	-40.1235	\\
		14	&	0.4022	&	0.5130	&	-59.9679	&	-42.9334	\\
		\hline
	\end{tabular}
	\centering
	\label{tab-3}
\end{table}

The simulation result implies that two PF solution curves have been obtained for the IEEE-14 bus system. However, the presented two continuous solutions will not be observed for the studied IEEE 14-bus system in practical if 1) reactive power generation in bus 8 exceeds its limit, i.e. $Q^{g,\text{max}}_8<6.7448$, 2) $P^d_7\not=0$ or $Q^d_7\not=0$ meaning $|V_7|=0$ is impossible to hold or 3) the voltage magnitude of any bus or power flow of any line exceeds its limit etc.

\section{Conclusion}
In this letter, two continuous PF solutions, or two PF solution curves, are obtained for the IEEE-14 bus system based on the SDP relaxation based method to locate all PF solutions even the slack bus is fixed. We admit that the observation is non-trivial. However, it is of great interest to researchers investigating the geometry or multiplicity nature of PF equations and more work, e.g. the conditions that PF solution curve appears, the relations between isolated and continuous PF solutions, may be of interest in future work.

\end{document}